%% file: aless_042516_final.tex
\documentclass{article}
\usepackage{amsfonts,amsmath,amssymb}
\usepackage{graphics,epsfig}
\usepackage{longtable,geometry}
\allowdisplaybreaks
\usepackage{color}
\allowdisplaybreaks
\makeatletter
\begin{document}
\input dibe_091511

\input wiener_mac

\title{Some Properties of DeGiorgi Classes}
\author{Emmanuele DiBenedetto\footnote{Supported by NSF grant DMS-1265548}\\
Department of Mathematics, Vanderbilt University\\  
1326 Stevenson Center, Nashville TN 37240, USA\\
email: {\tt em.diben@vanderbilt.edu}
\and
Ugo Gianazza\\
Dipartimento di Matematica ``F. Casorati", Universit\`a di Pavia\\ 
via Ferrata 1, 27100 Pavia, Italy\\
email: {\tt gianazza@imati.cnr.it}}
\maketitle
\begin{abstract}
The DeGiorgi classes $[DG]_p(E;\gm)$, defined in 
(1.1)${}_{\pm}$ below encompass, solutions of quasilinear 
elliptic equations with measurable coefficients as 
well as minima and Q-minima of variational integrals. For 
these classes we present some new results 
(\S~\ref{S:2} and \S~\ref{S:3:1}), and some known 
facts scattered in the literature 
(\S~\ref{S:3}--\S~\ref{S:5}), and formulate 
some open issues (\S~\ref{S:6}).
\vskip.2truecm
\noindent{\bf Mathematics Subject Classification (2010):} 
Primary 35J15, 49N60; Secondary 35J92
\vskip.2truecm
\noindent{\bf Key Words:} DeGiorgi classes, H\"older 
continuity, Harnack inequality, higher integrability, 
boundary behavior, decay estimates.
\end{abstract}
\section{Introduction}\label{S:1}
Let $E$ be open subset of $\rn$ and for $y\in\rn$, 
let $K_\rho(y)$ denote a cube of edge $2\rho$ 
centered at $y$. The DeGiorgi classes 
$[DG]_p^{\pm}(E;\gm)$ in $E$ are the 
collection of functions $u\in W^{1,p}_{loc}(E)$, 
for some $p>1$, satisfying
\begin{equation*}
\int_{K_\rho(y)}|D(u-k)_\pm|^pdx\le\frac{\gm}{(R-\rho)^p}
\int_{K_R(y)}|(u-k)_\pm|^pdx
\tag*{(1.1)${}_{\pm}$}
\end{equation*}
for all cubes $K_\rho(y)\subset K_R(y)\subset E$, and all 
$k\in\rr$, for a given positive constant $\gm$.
We further define
\begin{equation}\label{Eq:1:2}
[DG]_p(E;\gm)=[DG]_p^+(E;\gm)\cap[DG]_p^-(E;\gm).
\end{equation}
A celebrated theorem of DeGiorgi \cite{DeG} states that 
functions $u\in [DG]_p(E:\gm)$ are locally 
bounded and locally H\"older continuous in $E$. Moreover, 
non-negative functions $u\in[DG]_p(E:\gm)$ 
satisfy the Harnack inequality \cite{DbTru}. 

Local sub(super)-solutions, in $W^{1,p}_{\loc}(E)$, 
of quasi-linear elliptic equations in divergence form 
belong to $[DG]_p^{+(-)}(E;\gm)$ (\cite{LU}), 
with $\gm$ proportional to the ratio of upper and 
lower modulus of ellipticity. 
Local minima and/or $Q$-minima of variational integrals 
with $p$-growth with respect to $|Du|$ belong to these 
classes (\cite{GG}). Thus the $[DG]_p$-classes include 
local solutions of ellitic equations with merely 
bounded and measurable coefficients, only subject to some 
upper and lower ellipticity condition. They also 
include local minima or $Q$-minima of rather 
general functionals, even if not admitting a 
Euler equation. 

The interest in the DeGiorgi classes stems from 
the large class of, seemingly unrelated functions 
they encompass, and from properties, such as 
local H\"older continuity (\cite{DeG}), and the 
Harnack inequality (\cite{DbTru}), typically regarded 
as properties of solutions of elliptic partial 
differential equations (\cite{moser,LU}). 

The purpose of this note is to present some 
new results on DeGiorgi classes (\S~\ref{S:2} and \S~\ref{S:3:1}), 
as well as collecting some known facts scattered in the 
literature (\S~\ref{S:3}--\S~\ref{S:5}), and formulate 
some open issues (\S~\ref{S:6}) to serve as a basis 
for further investigations.
\section{DeGiorgi Classes and Sub(Super)-Harmonic 
Functions}\label{S:2}
\setcounter{equation}{1}
The generalized DeGiorgi classes $[GDG]_p^{\pm}(E;\gm)$, 
are the collection of functions $u\in W^{1,p}_{\loc}(E)$, 
for some $p>1$, satisfying 
\begin{equation*}
\int_{K_\rho(y)}|D(u-k)_\pm|^pdx\le\frac{\gm}{(R-\rho)^p}
\Big(\frac{R}{R-\rho}\Big)^{Np}\int_{K_R(y)}|(u-k)_\pm|^pdx
\tag*{(2.1)${}_{\pm}$}
\end{equation*}
for all cubes $K_\rho(y)\subset K_R(y)\subset E$, and all 
$k\in\rr$, for a given positive constant $\gm$.  Convex, 
monotone, non-decreasing functions of 
sub-harmonic functions are sub-harmonic. Similarly, 
concave, non-decreasing, functions of super-harmonic 
functions are super-harmonic. Similar statements hold 
for weak, sub(super)-solutions of linear ellipic equations 
with measurable coefficients (\cite{moser}).  
The next lemma establishes analogous properties 
for functions $u\in [DG]^{\pm}(E;\gm)$. 
Given any such class, we refer to the set of 
parameters $\{p,\gm,N\}$ as the {\it data} 
and say that a constant $C=C(\text{data})$ depends only 
on the {\it data} if it can be quantitatively determined 
a-priori only in terms of the indicated set of parameters.
\begin{lemma}\label{Lm:2:1} 
Let $\vp:\rr\to\rr$ be convex and non-decreasing, 
and let $u\in [DG]_p^+(E;\gm)$. 
There exists a positive constant $\ov{\gm}$ 
depending only on the data, and 
independent of $u$, such that 
$\vp(u)\in [GDG]_p^+(E;\ov{\gm})$. 

Likewise let $\psi:\rr\to\rr$ be concave and non-decreasing, 
and let $u\in [DG]_p^-(E;\gm)$. 
There exist a positive constant $\ov{\gm}$ 
depending only on the data, and 
independent of $u$, such that 
$\psi(u)\in [GDG]_p^-(E;\ov{\gm})$. 
\end{lemma}
\noi{\bf Proof:\ } By DeGiorgi's theorem (\cite{DeG,LU}), 
there exists a constant $C=C(\text{data})$, such that 
for any $u\in [DG]_p^{\pm}(E;\gm)$, there holds 
\begin{equation}\label{Eq:2:2}
\|(u-k)_{\pm}\|_{\infty, K_\rho(y)}\le\frac{C}{(R-\rho)^N}
\int_{K_R(y)}(u-k)_{\pm}dx
\end{equation}
for every pair of cubes $K_{\rho}(y)\subset K_R(y)\subset E$ 
and all $k\in\rr$. It suffices to prove the first statement 
for $\vp\in C^2(\rr)$, and verify that $\vp(u)$ 
satisfies {(2.1)${}_{+}$} for cubes $K_\rho\subset K_R$ 
centered at the origin of $\rn$. For any such $\vp$ and 
all $h\le k$
\begin{equation}\label{Eq:2:3}
\big(\vp(u)-\vp(h)\big)_+-\vp^\prime(h)(u-h)_+
=\int_{\rr^+}(u-k)_+\chi_{[k>h]} \vp^{\db}(k)dk
\end{equation}
From this, a.e. in $E$
\begin{equation*}
\big|D\big[\big(\vp(u)-\vp(h)\big)_+
-\vp^\prime(h) (u-h)_+\big]\big|^p
\le\Big(\int_{\rr}|D(u-k)_+|\chi_{[k>h]}\vp^{\db}(k) dk\Big)^p.
\end{equation*}
Integrate over $K_\rho$, take the $p$ root of both sides, 
and majorize the resulting term on the right-hand 
first by the continuous version of Minkowski 
inequality, then by applying the definition 
{(1.1)${}_{+}$} of the $[DG]_p^+(E;\gm)$-classes, and finally 
by using (\ref{Eq:2:2}). This gives
\begin{equation*}
\begin{aligned}
\big\|D\big[\big(\vp(u)-\vp(h)\big)_+
&-\vp^\prime(h) (u-h)_+\big]\big\|_{p,K_\rho}\\
&\le\int_{\rr}\|D(u-k)_+\|_{p,K_\rho}\hk\vp^{\db}(k) dk\\
&\le\frac{C}{R-\rho}\int_{\rr}
\|(u-k)_+\|_{p,K_{\frac{R+\rho}2}}\chi_{[k>h]}\vp^{\db}(k) dk\\
&\le\frac{C R^{\frac{N}{p}}}{R-\rho}\int_{\rr} 
\|(u-k)_+\|_{\infty, K_{\frac{R+\rho}2}} 
\chi_{[k>h]}\vp^{\db}(k) dk\\
&\le\frac{C R^{\frac{N}{p}}}{(R-\rho)^{N+1}}
\int_{\rr}\Big(\int_{K_R}(u-k)_+dx\Big)
\chi_{[k>h]}\vp^{\db}(k) dk\\
&=\frac{C R^{\frac{N}{p}}}{(R-\rho)^{N+1}}
\int_{K_R}\Big(\int_{\rr}(u-k)_+\chi_{[k>h]}\vp^{\db}(k) dk\Big)dx\\
&=\frac{C R^{\frac{N}{p}}}{(R-\rho)^{N+1}}
\int_{K_R}\big[\big(\vp(u)-\vp(h)\big)_+
-\vp^\prime(h) (u-h)_+\big]dx\\
&\le\frac{C}{R-\rho}\Big(\frac{R}{R-\rho}\Big)^N
\big\|\big(\vp(u)-\vp(h)\big)_+
-\vp^\prime(h) (u-h)_+\big\|_{p,K_R}.
\end{aligned}
\end{equation*}
In these calculations, we have denoted by $C=C(p,N,\gm)$ 
a generic constant depending only upon the data, and 
that might be different from line to line.  
In the last two steps we have interchanged the 
order of integration with the help of Fubini's Theorem and 
have applied H\"older's inequality.  By the convexity and 
monotonicity of $\vp$,
\begin{equation}\label{Eq:2:4}
\big(\vp(u)-\vp(h)\big)_+\ge\vp^\prime(h) (u-h)_+\ge0.
\end{equation}
Therefore, 
\begin{equation*}
\begin{aligned}
\big\|D\big(\vp(u)-\vp(h)\big)_+\big\|_{p,K_\rho}
&\le\frac{C}{R-\rho}\Big(\frac{R}{R-\rho}\Big)^N
\big\|\big(\vp(u)-\vp(h)\big)_+\big\|_{p,K_R}\\
&\quad +\big\|\vp^\prime(h) D(u-h)_+\big\|_{p,K_\rho}
\end{aligned}
\end{equation*}
Upon applying the definition of {(1.1)${}_{+}$} 
of $[DG]_p^+(E;\gm)$, and then (\ref{Eq:2:4}), 
the last term on the right-hand side is majorized by 
\begin{equation*}
\frac{C}{R-\rho}\big\|\vp(u)-\vp(h)_+\big\|_{p,K_R}.
\end{equation*}
Combining these estimates yields
\begin{equation}\label{Eq:2:5}
\int_{K_\rho(y)}\big|D\big(\vp(u)-k\big)_+\big|^pdx\le
\frac{\ov{\gm}}{(R-\rho)^p}\Big(\frac{R}{R-\rho}\Big)^N
\int_{K_R(y)}\big(\vp(u)-k\big)_+^pdx
\end{equation}
for all $k\in\rr$ and all $K_\rho(y)\subset K_R(y)\subset E$, 
for a constant $\ov{\gm}=\ov{\gm}(\text{data})$. \hfill\bbox
\vskip0.3cm

If $u\in [DG]_p^-(E;\gm)$ and $\vp$ is convex, 
there is no guarantee, in general, that 
$\vp(u)\in [GDG]_p^+(E;\ov{\gm})$ for some 
$\ov{\gm}=\ov{\gm}(p,N,\gm)$. The next lemma 
provides some sufficient conditions on $\vp$ 
for this to occur.  
\begin{lemma}\label{Lm:2:2}
Let $\vp:(a,+\infty)\to\rr$, for some $a<\infty$ be convex, 
non-increasing, and such that 
\begin{equation}\label{Eq:2:6}
\lim_{t\to+\infty}\vp(t)=\lim_{t\to+\infty}t\vp'(t)=0,
\end{equation}
and let $u\in [DG]_p^-(E;\gm)$, with range in $(a,+\infty)$. 
There exists a positive constant $\ov{\gm}$ depending only 
on the data, such that 
$\vp(u)\in [GDG]_p^+(E;\ov{\gm})$.  

Likewise let $\psi:(-\infty, a)\to\rr$, for some $a>-\infty$, 
be concave, non-increasing, and satisfying 
\begin{equation}\label{Eq:2:7}
\lim_{t\to-\infty}\psi(t)=\lim_{t\to-\infty}t\psi'(t)=0,
\end{equation}
and let $u\in [DG]_p^+(E;\gm)$, with range in $(-\infty,a)$. 
There exists a positive constant $\ov{\gm}$ depending only 
on the data, such that 
$\psi(u)\in [GDG]_p^-(E;\ov{\gm})$. 
\end{lemma}
\noi{\bf Proof:\ } It suffices to prove the first statement 
for $\vp\in C^2(\rr)$ over congruent cubes 
$K_\rho\subset K_R$ centered at the origin.  
The starting point is the analog of (\ref{Eq:2:3}), i.e., 
\begin{equation}\label{Eq:2:8}
\vp(u)=\int_{\rr}(u-k)_-\phi''(k)dk.
\end{equation}
Since $u\in [DG]_p^-(E;\gm)$, by (\ref{Eq:2:2}) the function $u$ 
is locally bounded below in $E$, and without loss of generality 
we may assume $u\ge0$. Hence the representation (\ref{Eq:2:8}) is 
well defined by virtue of the assumption (\ref{Eq:2:6}) on $\vp$. 
From this, by taking the gradient of both sides, then taking 
the $p$-power, and finally integrating over $K_\rho$ gives 
\begin{equation*}
\int_{K_{\rho}}|D\phi(u)|^p dx=\int_{K_{\rho}}\Big|
\int_{\rr^+}D(u-k)_-\phi''(k)dk\Big|^pdx.
\end{equation*}
The proof now parallels that of Lemma~\ref{Lm:2:1}. 
Specifically, apply sequentially the continuous version of 
Minkowski's inequality, the definition {(1.1)${}_-$} 
of the classes $[DG]_p^-(E;\gm)$, the sup-bound 
(\ref{Eq:2:2}), interchange the order of integration, 
and use H\"older's inequality. This gives 
\begin{align*}
\|D\vp(u)\|_{p,K_\rho}&\le\int_{\rr^+}\|D(u-k)_-\|_{p,K_\rho}
\vp^{\db}(k) dk\\
&\le\frac{C}{R-\rho}\int_{\rr^+}\|(u-k)_-\|_{p,K_{\frac{R+\rho}2}}
\vp^{\db}(k)dk\\
&\le\frac{C R^{\frac{N}p}}{R-\rho}
\int_{\rr^+}\|(u-k)_-\|_{\infty,K_{\frac{R+\rho}2}}\vp^{\db}(k)dk\\
&\le\frac{C R^{\frac{N}p}}{(R-\rho)^{N+1}}
\int_{\rr^+}\int_{K_R}(u-k)_-\vp^{\db}(k)dk\\
&=\frac{C R^{\frac{N}p}}{(R-\rho)^{N+1}}
\int_{K_R}\vp(u)dx\\
&=\frac{C}{(R-\rho)}\Big(\frac{R}{R-\rho}\Big)^N\|\vp(u)\|_{p, K_R}.
\end{align*}
Now if $\vp$ is convex, non-increasing and satisfying 
(\ref{Eq:2:6}), the function $(\vp-\ell)_+$, for all $\ell$ 
in the range of $\vp$, shares the same properties. Hence 
\begin{equation*}
\int_{K_{\rho}(y)} \big|D\big(\vp(u)-\ell\big)_+\big|^{p}dx \le 
\frac{C}{(R-\rho)^p}\Big(\frac{R}{R-\rho}\Big)^{Np}
\int_{K_R(y)}\big(\vp(u)-\ell\big)_+^{p}dx
\end{equation*}
for all cubes $K_{\rho}(y)\subset K_R(y)\subset E$ and 
all $\ell \in\rr$.\hfill\bbox
\subsection{Some Consequences}\label{S:2:1}
The sup-bound in (\ref{Eq:2:2}) can be given the following 
sharper form (\cite{DbTru}). 
\begin{lemma}\label{Lm:2:3}
Let $u\in [DG]_p^{\pm}(E;\gm)$. Then for all $\sig>0$ 
there exists a constant $C_\sig$ depending only upon 
the data and $\sig$, such that 
\begin{equation}\label{Eq:2:9}
\sup_{K_\rho(y)}(u-k)_{\pm}\le C_\sig
\Big(\frac{R}{R-\rho}\Big)^{\frac{N}{\sig}} 
\Big(\dashint_{K_R(y)}(u-k)_{\pm}^\sig dx\Big)^{\frac1{\sig}}.
\end{equation}
\end{lemma}
If $u\in[DG]_p^-(E;\gm)$ is non-negative, then 
Lemma~\ref{Lm:2:2} with $\vp(u)=u^{-1}$ and $a=0$, 
implies that $u^{-1}\in [GDG]_p^+(E;\gm)$.  Therefore 
Lemma~\ref{Lm:2:3}, with $k=0$, implies that for 
all $\tau>0$, 
\begin{equation}\label{Eq:2:10}
\frac1{\dsty\inf_{K_\rho(y)}u}\le C_\tau
\Big(\frac{R}{R-\rho}\Big)^{\frac{N}{\tau}}
\Big(\dashint_{K_R(y)}\frac1{u^\tau} dx\Big)^{\frac1{\tau}}.
\end{equation}
\begin{proposition}\label{Prop:2:1}
Let $u$ be a non-negative function in the DeGiorgi 
classes $[DG]_p(E;\gm)$. Then for any pair of positive
numbers $\sig$ and $\tau$
\begin{equation}\label{Eq:2:11}
\frac{\dsty\sup_{K_\rho(y)}u}{\dsty\inf_{K_\rho(y)}u}\le C_{\sig}C_{\tau}
\Big(\frac{R}{R-\rho}\Big)^{N(\frac1{\sig}+\frac1{\tau})} 
\Big(\dashint_{K_R(y)}u^\sig dx\Big)^{\frac1{\sig}}
\Big(\dashint_{K_R(y)}\frac1{u^\tau} dx\Big)^{\frac1{\tau}}.
\end{equation}
\end{proposition}
Inequalities of the form (\ref{Eq:2:9}) are at the basis of 
Moser's approach to the Harnack inequality for non-negative 
weak solutions to quasilinear elliptic equations with bounded 
and measurable coefficients (\cite{moser}). The Harnack 
inequality will follow from (\ref{Eq:2:11}) if $\ln u\in BMO(E)$. 
This fact is established by Moser for non-negative 
weak solutions of elliptic equations. We will establish that 
for non-negative functions $u\in [DG]_p^-(E;\gm)$, one has 
$\ln u\in BMO(E)$ by using the Harnack inequality 
established in (\cite{DbTru}). 
\section{DeGiorgi Classes, $BMO(E)$ and Logarithmic 
Estimates}\label{S:3}
The proof of the following lemma is in \cite{DbTru}.
\begin{lemma}\label{Lm:3:1}
Let $u\in[DG]_p^-(E;\gm)$ be non-negative. There exist 
positive constants $C$ and $\sig$, depending only upon 
the data, such that 
\begin{equation}\label{Eq:3:1}
\dashint_{K_\rho(y)} u^{\sig}dx\le C\inf_{K_\rho(y)} u^{\sig},
\end{equation}
for any pair of cubes $K_\rho(y)\subset K_{2\rho}(y)\subset E$. 
\end{lemma}
Such an inequality, referred to as the weak Harnack inequality, was 
established by Moser for non-negative super-solutions of 
elliptic equations with bounded and measurable coefficients 
(\cite{moser}). It is noteworthy that it continues to hold 
for non-negative functions in $[DG]_p^-(E;\gm)$, with no further 
reference to equations. 
\begin{lemma}\label{Lm:3:2}
Let $u\in[DG]_p^-(E;\gm)$ be non-negative. Then $\ln u\in BMO$.
\end{lemma}
\ms
\noi{\bf Proof:\ } By Lemma~\ref{Lm:3:1}
\begin{equation}\label{Eq:3:2}
\begin{aligned}
\dashint_{K_\rho(y)} u^\sig dx\,\dashint_{K_\rho(y)}\frac1{u^\sig}dx 
&\le \dashint_{K_\rho(y)} u^\sig dx\,\sup_{K_\rho(y)}\frac1{u^{\sig}}\\
&=\dashint_{K_\rho(y)} u^\sig dx\,\frac1{\dsty\inf_{K_\rho(y)} u^\sig}\le C
\end{aligned}
\end{equation}
for any pair of cubes $K_\rho(y)\subset K_{2\rho}(y)\subset E$. Set
\begin{equation*}
(\ln u^{\sig})_\rho=\dashint_{K_\rho(y)}\ln u^{\sig} dx,
\end{equation*}
and estimate
\begin{equation*}
\begin{aligned}
\dashint_{K_\rho(y)}e^{|\ln u^\sig-(\ln u^\sig)_\rho|}dx&\le 
e^{-(\ln u^\sig)_\rho}\>\dashint_{K_\rho(y)} e^{\ln u^\sig}dx\\
&\quad +e^{(\ln u^\sig)_\rho}\>\dashint_{K_\rho(y)} e^{-\ln u^\sig}dx.
\end{aligned}
\end{equation*}
The second term on the right-hand side is estimated by Jensen's 
inequality and (\ref{Eq:3:2}) and yields
\begin{equation*}
\begin{aligned}
e^{(\ln u^\sig)_\rho}\>\dashint_{K_\rho(y)} e^{-\ln u^\sig}dx&\le 
\dashint_{K_\rho(y)}e^{\ln u^\sig}dx\>
\dashint_{K_\rho(y)}\frac1{u^\sig} dx\\
&\le \dashint_{K_\rho(y)} u^\sig dx\>
\dashint_{K_\rho(y)} \frac1{u^\sig} dx\le C 
\end{aligned}
\end{equation*}
The first term is estimated analogously. Hence, there exists 
a constant $\bar{C}$, depending only upon the data, such that
\begin{equation*}
\dashint_{K_\rho(y)}e^{|\ln u^\sig-(\ln u^\sig)_\rho|}dx\le \bar{C} 
\end{equation*}
for any pair of cubes $K_\rho(y)\subset K_{2\rho}(y)\subset E$. Thus 
$\ln u\in BMO(E)$.\hfill\bbox
\subsection{Logarithmic Estimates Revisited}\label{S:3:1}
Let $u\in W_{\loc}^{1,p}(E)$ be a non-negative weak super-solution 
of an elliptic equation in divergence form, and with only bounded 
and measurable coefficients. Then there exists a constant $C$, depending 
only on $p$, $N$, and the modulus of ellipticity of the equation, 
such that 
\begin{equation}\label{Eq:3:3}
\dashint_{K_\rho(y)}|D\ln u|^pdx\le \frac{C}{(R-\rho)^p}
\end{equation}
for every pair of cubes $K_\rho(y)\subset K_R(y)\subset E$. Such 
an estimate, established by Moser, permits one to prove that 
$\ln u\in BMO(E)$, which in turn yields the Harnack inequality. Our 
approach for functions in the $[DG]_p^-(E;\gm)$ classes is somewhat 
different. For non-negative functions in such classes we first 
establish the weak Harnack estimate (\ref{Eq:3:1}), and then the 
latter is used to prove Lemma~\ref{Lm:3:2}. It is not known, whether 
non-negative functions in $[DG]_p^-(E;\gm)$ satisfy (\ref{Eq:3:3}). The 
next proposition is a partial result in this direction. 
\begin{proposition}\label{Prop:3:1}
Let $u\in [DG]_p^-(E;\gm)$ be non-negative and bounded above 
by some positive constant $M$. Then
\begin{equation}\label{Eq:3:4}
\int_{K_\rho(y)}|D\ln u|^pdx\le \frac{\gm p}{(R-\rho)^p}\int_{K_R(y)} 
\ln{\frac{M}{u}}dx
\end{equation}
for any pair of cubes $K_\rho(y)\subset K_{R}(y)\subset E$. 
\end{proposition}
\ms
\noi{\bf Proof:\ } The arguments being local may assume 
that $y=\{0\}$. By the definition {(1.1)${}_-$} classes, for all 
$0<t<M$, 
\begin{equation*}
\int_{K_\rho}|D(u-t)_-|^p dx\le\frac{\gm}{(R-\rho)^p}
\int_{K_R}(u-t)_-^p dx.
\end{equation*}
Multiply both sides by $t^{-p-1}$ and integrate over $(0,M)$. The 
 left-hand side is transformed as
\begin{align*}
\int_0^M{\frac{dt}{t^{p+1}}}\int_{K_\rho}|D(u-t)_-|^pdx&=
\int_{K_\rho}\Big(\int_0^M|D(u-t)_-|^p{\frac1{t^{p+1}}}dt\Big)dx\\
&=\int_{K_\rho}|Du|^p\Big(\int_0^M{\frac1{t^{p+1}}}
\chi_{[u<t]}dt\Big)dx\\
&=\int_{K_\rho}|Du|^p\Big(\int_u^M{\frac1{t^{p+1}}}dt\Big)dx\\
&=\int_{K_\rho}\Big(-{\frac1p}{\frac{|Du|^p}{M^p}}
+{\frac1p}{\frac{|Du|^p}{u^p}}\Big)dx\\
&={\frac1p}\int_{K_\rho}|D\ln u|^pdx-{\frac1{pM^p}}
\int_{K_\rho}|Du|^pdx.
\end{align*}
The integral on the right-hand side is transformed as
\begin{align*}
\int_0^M\frac1{t^{p+1}}\Big(\int_{K_R}(u-t)_-^pdx\Big)dt&=
\int_{K_R}\Big(\int_u^M\frac{(t-u)^p}{t^{p+1}}dt\Big)dx\\
&=\int_{K_R}\Big[-\frac1p\frac{(t-u)^p}{t^p}\Big|_u^M
+\int_u^M \frac{(t-u)^{p-1}}{t^{p-1}}\frac{dt}t\Big]dx\\
&=-\frac1{pM^p}\int_{K_R}(M-u)^pdx+
\int_{K_R}\Big(\int_u^M \Big(\frac{t-u}t\Big)^{p-1}\frac{dt}t\Big)dx\\
&\le -\frac1{pM^p}\int_{K_R}(M-u)^pdx+\int_{K_R}
\ln{\frac Mu}dx.
\end{align*}
Combining the previous estimates gives
\begin{align*}
\int_{K_\rho}|D\ln u|^p dx&\le{\frac1{M^p}}\Big(
\int_{K_\rho}|Du|^p dx-\frac{\gm}{(R-\rho)^p}
\int_{K_R}(M-u)^pdx\Big)\\
&\quad+\frac{\gm p}{(R-\rho)^p}\int_{K_R}\ln{\frac Mu}dx.
\end{align*}
Since $u\in [DG]_p^-(E;\gm)$, the term in round brackets on 
the right-hand side is non-positive and can be discarded.\hfill\bbox
\begin{remark}\label{Rmk:3:1} {\normalfont
Applying Lemma~\ref{Lm:2:2} to $\dsty\phi(u)=\ln_+(M/u)$, gives 
the weaker estimate
\begin{equation}\label{Eq:3:5}
\int_{K_\rho(y)}|D\ln u|^p dx\le\frac{\ov{\gm}}{(R-\rho)^p}
\int_{K_R(y)}\Big(\ln{\frac Mu}\Big)^p dx.
\end{equation}
}
\end{remark}
\section{Higher Integrability of the Gradient of Functions 
in the DeGiorgi Classes}\label{S:4}
\begin{proposition}\label{Prop:4:1}
Let $u\in\dg$. Then there exist constants $C>1$ and $\sig>0$, 
dependent only upon the data, such that, for any pair of cubes 
$K_\rho(y)\subset K_{R}(y)\subset E$, there holds
\begin{equation}\label{Eq:6:5}
\Big(\dashint_{K_{\rho}(y)}|Du|^{p(1+\sig)}dx\Big)^{\frac1{p(1+\sig)}}
\le C\Big(\frac{R}{\rho}\Big)^{\frac{N}{p}}\Big(\frac{R}{R-\rho}\Big)
\Big(\dashint_{K_R(y)}|Du|^pdx\Big)^{\frac1p}.
\end{equation}
\end{proposition}
\noi{\bf Proof:\ } Let $u$ be in the classes $[DG]_p(E;\gm)$ 
defined in (\ref{Eq:1:2}). For any pair of cubes 
$K_\rho(y)\subset K_R(y)\subset E$, write 
down {(1.1)${}_+$} and {(1.1)${}_-$} for the choice 
\begin{equation*}
k=u_R\df{=}\dashint_{K_R(y)}u dx.
\end{equation*}
Adding the resulting inequalities gives
\begin{equation*}
\int_{K_\rho(y)}|Du|^pdx\le\frac{\gm}{(R-\rho)^p}\int_{K_R(y)}|u-u_R|^pdx.
\end{equation*}
By the Sobolev-Poincar\'e inequality
\begin{equation*}
\dashint_{K_R(y)}|u-u_R|^pdx\le C_q\,R^p
\Big(\dashint_{K_R(y)}|Du|^qdx\Big)^{\frac pq},\qquad\text{ for all }\> 
q\in\Big[\frac{Np}{N+p},p\Big]
\end{equation*}
for a constant $C_q=C_q(N,q)$. Hence for all such $q$
\begin{equation*}
\dashint_{K_\rho(y)}|Du|^pdx\le C_q \gm
\Big(\frac{R}{R-\rho}\Big)^p\Big(\frac{R}{\rho}\Big)^N
\Big(\dashint_{K_R(y)}|Du|^qdx\Big)^{\frac pq}
\end{equation*}
for all pair of congruent cubes $K_\rho(y)\subset K_R(y)\subset E$. 
The conclusion follows from this and the local version of 
Gehring's lemma (\cite{gehring}), as appearing in \cite{giamo}.\hfill\bbox
\begin{remark}\label{Rmk:4:1} {\normalfont
Hence the higher integrability of the gradient of solutions 
of elliptic equations with measurable coefficients (\cite{stredulinsky}), 
and more generally of $Q$-minima (\cite{GG}), continues to 
hold for function in the DeGiorgi classes.
If $u\in [DG]_p^{\pm}(E;\gm)$, the conclusion 
is in general false, as one can verify starting 
from sub(super)-harmonic functions. However, 
essentially the same arguments give the inequality 
\begin{gather*}
\dashint_{K_\rho(y)}|D(u-k)_{\pm}|^pdx\le C_q \gm
\Big(\frac{R}{R-\rho}\Big)^p\Big(\frac{R}{\rho}\Big)^N
\Big(\dashint_{K_R}|Du|^qdx\Big)^{\frac pq}\quad\text{for all }\> q\in\Big[\frac{Np}{N+p},p\Big],\\ 
\text{ and all }\>
k\ge\dashint_{K_R(y)} udx\ \  \text{ if }\ u\in [DG]_p^{+}(E;\gm),\quad
k\le\dashint_{K_R(y)} udx\ \ \text{ if }\ u\in [DG]_p^{-}(E;\gm).
\end{gather*}
}
\end{remark} 
\section{Measure Theoretical Decay Estimates of Functions 
in DeGiorgi Classes}\label{S:5}
For a non-negative function $f\in L^1_{\loc}(E)$ one estimates 
the measure of the set $[f>t]$ relative to a cube 
$K_\rho(y)\subset E$, as $\mu\big([f>t]\cap K_\rho(y)\big)
\le t^{-1}\|f\|_{1,K_\rho(y)}$. Estimates of the measure of 
the set $[f<t]$ relative to $K_\rho(y)$ are not, in general, 
a consequence of the mere integrability of $f$. One of DeGiorgi's 
estimates of \cite{DeG}, is that if $u$ is a non-negative function 
in $[DG]_p^-(E;\gm)$, then 
\begin{equation}\label{Eq:5:1}
\frac{\big|[u<t]\cap K_\rho(y)\big|}{|K_\rho|}\le
\frac{C(N,p,\gm)}{|\ln t|^{1/p}} 
\qquad\text{ asymptotically as }\>t\to0,
\end{equation}
provided $|[u>t]\cap K_\rho(y)|\ge\frac12|K_\rho|$.
Here $|\Sig|$ denotes the Lebesgue measure of a measurable 
set $\Sig\subset\rn$. The next proposition improves on this estimate.
\begin{proposition}\label{Prop:5:1} 
Let $u\in [DG]_p^-(E;\gm)$ be non-negative, and assume that for some 
$t_o>0$ and $\al\in(0,1)$, there holds 
\begin{equation}\label{Eq:5:2}
\frac{\big|[u>t_o]\cap K_\rho(y)]\big|}{|K_\rho|}\ge\al.
\end{equation}
There exist positive constants $C,t_*,\sig=C,t_*,\sig(N,p,\gm,t_o,\al)$, 
depending only on the indicated parameters and independent of $u$, 
such that 
\begin{equation}\label{Eq:5:3}
\frac{\big|[u<t]\cap K_\rho(y)\big|}{|K_\rho|}\le
\frac{C}{|\ln t|^{\sig|\ln t|^{\frac12}}}, \qquad\text{ for }\>t<t_*.
\end{equation}
\end{proposition}
\vskip0.3cm
\noi{\bf Proof:\ } 
In what follows we denote by $C$ a generic positive 
constant that can be determined a-priori only in terms of 
$\{N,p,\gm,t_o,\al\}$ and that it may be different 
in the same context. The arguments being local to 
concentric cubes $K_{\rho}(y)\subset K_{2\rho}(y)\subset E$, 
may assume $y=\{0\}$ and write $K_\rho(0)=K_\rho$. 
Let $n_o$ be the smallest positive integer such 
that $2^{-n_o}\le t_o$, and for $n\ge n_o$ set 
\begin{equation*}
A_{n,\rho}\df{=}\Big[u<\frac1{2^n}\Big]\cap K_\rho,
\qquad\text{ for }\> n\ge n_o. 
\end{equation*}
The discrete isoperimetric inequality 
(\cite[Chapter~I, Lemma~2.2]{dibe-sv}), reads
\begin{equation*}
(\ell-h)\big|[u<h]\cap K_\rho\big|\le C(N)
\frac{\rho^{N+1}}{\big|[u>\ell]\cap K_\rho\big|}
\int_{[h<u<\ell]\cap K_\rho} |Du|dx
\end{equation*}
for any two levels $0<h<\ell$. Applying it with 
\begin{equation*}
\ell=\frac1{2^n},\quad h=\frac1{2^{n+1}},\quad\text{ so that }\quad 
[h<u<\ell]\cap K_\rho= A_{n,\rho}-A_{n+1,\rho},
\end{equation*}
and taking into account (\ref{Eq:5:2}), yields
\begin{equation*}
\frac1{2^{n+1}}\big|A_{n+1,\rho}\big|\le\frac{C(N)}{\al} \rho^N
\int_{A_{n,\rho}-A_{n+1,\rho}} |Du|dx.
\end{equation*}
Majorize the right-hand side by the H\"older inequality,  then 
raise both terms to the power $\frac p{p-1}$, and majorize 
the right-hand side by {(1.1)${}_-$} in the definition 
of the classes $[DG]_p^-(E;\gm)$. These sequential estimates yield
\begin{align*}
\frac{1}{2^{n\frac{p}{p-1}}}|A_{n+1,\rho}|^{\frac p{p-1}}
&\le C\rho^{\frac p{p-1}}\Big(\int_{K_\rho}
\big|D\big(u-{\txty\frac1{2^{n}}}\big)_-|^pdx\Big)^{\frac1{p-1}}
\big|A_{n,\rho}- A_{n+1,\rho}|\\
&\le C\Big(\int_{K_\rho}
\big(u-{\txty\frac1{2^{n}}}\big)_-^pdx\Big)^{\frac1{p-1}}
\big|A_{n,\rho}- A_{n+1,\rho}|\\
&\le \frac{C}{2^{n\frac{p}{p-1}}}|A_{n_o,2\rho}|^{\frac1{p-1}}
\big|A_{n,\rho}- A_{n+1,\rho}|.
\end{align*}
This in turn yields the recursive inequalities
\begin{equation*}
\big|A_{n+1,\rho}\big|^{\frac{p}{p-1}}
\le C(N,p,\gm,\al) \big|A_{n_o,2\rho}\big|^{\frac1{p-1}}
\big|A_{n,\rho}- A_{n+1,\rho}|.
\end{equation*}
Let $n_*$ be a positive integer to be chosen. Adding them from 
$n_o$ to $n_*-1$ gives
\begin{equation}\label{Eq:5:4}
\big|A_{n_*,\rho}\big|\le\frac{C(N,p,\gm,\al)}{(n_*-n_o)^{\frac{p-1}p}}
\big|A_{n_o,2\rho}\big|^{\frac1p}\big|A_{n_o,\rho}|^{\frac{p-1}p}.
\end{equation}
Return now to the assumption (\ref{Eq:5:2}) and estimate
\begin{equation*}
\frac{\big|[u>t_o]\cap K_{2\rho}(y)]\big|}{|K_{2\rho}|}\ge
\frac{\big|[u>t_o]\cap K_\rho(y)]\big|}{2^N|K_\rho|}\ge\frac{\al}{2^N}.
\end{equation*}
Therefore, the same arguments leading to (\ref{Eq:5:4}) can be repeated 
over the cube $K_{2\rho}$ and give
\begin{equation}\label{Eq:5:5}
\big|A_{n_*,2\rho}\big|\le\frac{C(N,p,\gm,\al)}{(n_*-n_o)^{\frac{p-1}p}}
\big|A_{n_o,4\rho}\big|^{\frac1p}\big|A_{n_o,2\rho}|^{\frac{p-1}p}.
\end{equation}
While the constant $C$ in (\ref{Eq:5:5}) differs from the one 
in (\ref{Eq:5:4}), we may take them to be equal by taking the 
largest. The assumption (\ref{Eq:5:2}) continue to hold with 
$t_o$ replaced by $2^{-n_*}$. Hence the previous arguments can be 
repeated and yield the analogues of (\ref{Eq:5:4})--(\ref{Eq:5:5}), i.e.,
\begin{equation*}
\begin{aligned}
\big|A_{2n_*,\rho}\big|&\le{\dsty\frac{C(N,p,\gm,\al)}{(n_*-n_o)^{\frac{p-1}p}}
\big|A_{n_*,2\rho}\big|^{\frac1p}\big|A_{n_*,\rho}|^{\frac{p-1}p}}\\
\big|A_{2n_*,2\rho}\big|&\le{\dsty\frac{C(N,p,\gm,\al)}{(n_*-n_o)^{\frac{p-1}p}}
\big|A_{n_*,4\rho}\big|^{\frac1p}\big|A_{n_*,2\rho}|^{\frac{p-1}p}}
\end{aligned}
\end{equation*}
for the same constant $C$. Combining them gives 
\begin{equation*}
\big|A_{2n_*,\rho}\big|\le\frac{C^2 4^{2N}}{(n_*-n_o)^{2\frac{p-1}p}}
|K_\rho|.
\end{equation*}
Iteration of this procedure yields
\begin{equation*}
\big|A_{jn_*,\rho}\big|\le\frac{C^j 4^{jN}}{(n_*-n_o)^{j\frac{p-1}p}}
|K_\rho|\qquad\text{ for all }\>j\in\nn.
\end{equation*}
Choose $n_*$ so large that $n_*-n_o>\frac12n_*$, and then take 
$j=n_*$. By possibly modifying the various constants, the previous 
inequality yields
\begin{equation*}
\big|A_{j^2,\rho}\big|\le\frac{C^j 4^{jN}}{j^{j\frac{p-1}p}}
|K_\rho|\qquad\text{ for all }\>j\in\nn.
\end{equation*}
The constant $C$ being fixed, for each $0<\varep<\frac{p-1}p$ there 
exists $j^*$ so large that
\begin{equation*}
\big|A_{j^2,\rho}\big|\le\frac{1}{j^{j\varep}}
|K_\rho|\qquad\text{ for all }\>j\ge j^*.
\end{equation*}
Fix now $t\le2^{-j^{*2}}$ and let $j$ be the largest integer 
such that $2^{-(j+1)^2}\le t\le 2^{-j^2}$. For such choices
\begin{equation*}
\frac{\big|[u<t]\cap K_\rho\big|}{|K_\rho|}\le 
\frac{\big|A_{j^2,\rho}\big|}{|K_\rho|}\le
\frac{C}{|\ln t|^{\frac{\varep}2|\ln t|^{\frac12}}}.
\tag*{\bbox}
\end{equation*}
The parabolic version of this result has been used in \cite{DBGV-indiana}.
\section{Boundary Behavior of Functions in the DeGiorgi Classes}\label{S:6}
Let $h\in W_{\loc}^{1,p}(\rn)\cap C(\rn)$. The DeGiorgi 
classes $[DG]_p^{+(-)}(\bar{E};\gm,h)$, in the closure of 
$E$ are the collection of functions $u\in W^{1,p}_{\loc}(\bar{E})$, 
such that $(u-h)\in W^{1,p}_o(E\cap K_R(y))$, for all cubes 
$K_R(y)$ centered at some $y\in\ple$, and satisfying 
\begin{equation}\label{Eq:6:1}
\int_{K_\rho(y)\cap E}|D(u-k)_{+(-)}|^pdx\le\frac{\gm}{(R-\rho)^p}
\int_{K_R(y)\cap E}(u-k)_{+(-)}^pdx
\end{equation}
for all pairs of congruent cubes $K_\rho(y)\subset K_R(y)$, centered 
at some $y\in\ple$ and all levels
\begin{equation}\label{Eq:6:2}
k\ge\sup_{K_R(y)\cap\ple}h,\qquad \Big( k\le\inf_{K_R(y)\cap\ple}h\Big).
\end{equation}
We let further
\begin{equation*}
[DG]_p(\bar{E};\gm,h)= [DG]_p^+(\bar{E};\gm,h)\cap [DG]_p^-(\bar{E};\gm,h).
\end{equation*}
Functions in $[DG]_p(\bar{E};\gm,h)$ are continuous up to 
points $y\in\ple$, provided $E$ satisfies a positive geometric 
density at $y$, i.e., there exist $\rho_o$ and $\eta\in(0,1)$, 
such that (see \cite{LU})
\begin{equation*}
\big|E^c\cap K_\rho(y)\big|\ge\eta|K_\rho(y)|,
\qquad\text{ for all }\> \rho\le\rho_o.
\end{equation*}
For $1<p<N$, the $p$-capacity of the 
compact set $E^c\cap\bar{K}_\rho(y)$ is defined by
\begin{equation}\label{Eq:6:3}
c_p[E^c\cap\bar{K}_\rho(y)]=\inf_{
\ttop{\psi\in W^{1,p}_o(\rn)\cap C(\rn)}{E^c\cap\bar{K}_\rho(y)
\subset[\psi\ge1]}}\int_{\rn}|D\psi|^pdx.
\end{equation}
For {$1<p< N$}, the relative $p$-capacity of 
$E^c\cap\bar{K}_\rho(y)$ with respect to $K_\rho(y)$ is
\begin{equation}\label{Eq:6:4}
\qquad 
\dl_y(\rho)=\frac{c_p[E^c\cap\bar{K}_\rho(y)]}{\rho^{N-p}}, 
\qquad\qquad (1<p<N).
\end{equation}
If $p=N$, and for $0<\rho<1$, the $N$-capacity of the 
compact set $E^c\cap\bar{K}_\rho(y)$, with respect to the 
cube $K_{2\rho}(y)$, is defined by 
\begin{equation}\label{Eq:6:6}
c_N[E^c\cap\bar{K}_\rho(y)]={\inf_{
\ttop{\psi\in W^{1,N}_o(K_{2\rho}(y))\cap 
C_o(K_{2\rho}(y))}{E^c\cap\bar{K}_\rho(y)
\subset[\psi\ge1]}}\int_{K_{2\rho}(y)}|D\psi|^Ndx}.
\end{equation}
The relative capacity $\dl_y(\rho)$ can be formally 
defined by (\ref{Eq:6:4}), for all $1<p\le N$. For $p=N$, we let 
$\dl_y(\rho)\equiv c_N[E^c\cap\bar{K}_\rho(y)]$, as 
defined by (\ref{Eq:6:6}).
For a positive parameter $\eps$ denote by $I_{p,\eps}(y,\rho)$ 
the Wiener integral of $\ple$ at $y\in\ple$, i.e.,   
\begin{equation}\label{Eq:6:7}
I_{p,\eps}(y,\rho)=\int_\rho^1[\dl_y(t)]^{\frac{1}{\eps}}
\frac{dt}{t}.
\end{equation}
The celebrated Wiener criterion states that a harmonic 
function in $E$ is continuous up to $y\in\ple$ if and only 
if the Wiener integral $I_{2,1}(y,\rho)$ diverges as 
$\rho\to0$ (\cite{wiener}). 

It is known that weak solutions of quasilinear equations 
in divergence form, and with principal part exhibiting a 
$p$-growth with respect to $|Du|$, when given continuous 
boundary data $h$ on $\ple$, are continuous up to $y\in\ple$ 
if $I_{p,(p-1)}(y,\rho)$ diverges as $\rho\to0$ (\cite{gar-zie}). 
Since such solutions belong to the boundary 
$[DG]_p(\bar{E};\gm,h)$ classes (\cite{GG}), it is natural 
to ask whether the divergence of the Wiener integral 
$I_{p,(p-1)}(y,\rho)$, is sufficient to insure the boundary 
continuity for functions $u\in[DG]_p(\bar{E};\gm,h)$. 

The only result we are aware of in this direction is due to 
Ziemer (\cite{ziem}). It states that a function 
$u\in[DG]_p(\bar{E};\gm,h)$ is continuous up to $y\in\ple$ if 
\begin{equation}\label{Eq:1:12}
\int_{\rho}^1 \exp\Big(-\frac1{\dl_y(t)^{\frac1{p-1}}}\Big)
\frac{dt}{t}\>\to\>\infty\quad\text{ as }\> \rho\to0.
\end{equation}
Ziemer's proof follows from a standard DeGiorgi iteration 
technique. It has been recently established that local 
minima of variational integrals when given 
continuous boundary data $h$ are continuous up to $y\in\ple$ 
provided (\cite{DiGi}) $I_{p,\varep}(y,\rho)$ diverges as 
$\rho\to0$. Here $\varep$ is a number that can be determined 
a-priori only in terms of the growth properties 
of the functional. While such minima are in the 
classes $[DG]_p(\bar{E};\gm,h)$, the result is not known to hold 
for functions merely in such classes. Also the optimal parameter 
$e=(p-1)$ remains elusive. A similar result has been recently 
obtained with a different approach in \cite{jana}.

The significance of a Wiener condition for Q-minima, 
is that the structure of $\ple$ near a boundary point 
$y\in\ple$, for $u$ to be continuous up to $y$, hinges 
on minimizing a functional, rather than solving 
an elliptic p.d.e. 
 
\bye

%% file: dibe_091511.tex
\newtheorem{theorem}{Theorem}[section]
\newtheorem{proposition}{Proposition}[section]
\newtheorem{lemma}{Lemma}[section]
\newtheorem{corollary}{Corollary}[section]
\newtheorem{remark}{Remark}[section]
\newtheorem{definition}{Definition}[section]
\renewcommand{\thesection}{\arabic{section}}
\renewcommand{\theequation}{\thesection.\arabic{equation}}
\renewcommand{\thetheorem}{\thesection.\arabic{theorem}}
\numberwithin{equation}{section}
\numberwithin{theorem}{section}
\numberwithin{proposition}{section}
\numberwithin{lemma}{section}
\numberwithin{remark}{section}
\setcounter{secnumdepth}{3}
\newcommand{\cl}{\centerline}
\newcommand{\sms}{\smallskip}
\newcommand{\ms}{\medskip}
\newcommand{\bs}{\bigskip}
\newcommand{\noi}{\noindent}
\newcommand{\itl}[1]{\textit{#1}}
\newcommand{\blf}[1]{\textbf{#1}}
\newcommand{\dsty}{\displaystyle}
\newcommand{\txty}{\textstyle}
\newcommand{\ssty}{\scriptstyle}
\newcommand{\tty}{\texttt}


\newcommand\Par{\mathhexbox278\,}


\newcommand{\al}{\alpha}
\newcommand{\Al}{\Alpha}
\newcommand{\be}{\beta}
\newcommand{\Be}{\Beta}
\newcommand{\Gm}{\Gamma}
\newcommand{\gm}{\gamma}
\newcommand{\dl}{\delta}
\newcommand{\Dl}{\Delta}
\newcommand{\lm}{\lambda}
\newcommand{\Lm}{\Lambda}
\newcommand{\kp}{\kappa}
\newcommand{\varep}{\varepsilon}
\newcommand{\eps}{\epsilon}
\newcommand{\vp}{\varphi}
\newcommand{\sig}{\sigma}
\newcommand{\Sig}{\Sigma}
\newcommand{\om}{\omega}
\newcommand{\Om}{\Omega}
\newcommand{\uom}{\mbox{\boldmath$\omega$}}
\newcommand{\btau}{\mbox{\boldmath$\tau$}}
\newcommand{\bnu}{\mbox{\boldmath$\nu$}}
\newcommand{\up}{\upsilon}
\newcommand{\z}{\zeta}


\newcommand{\df}[1]{\buildrel\mbox{\small def}\over{#1}}
\newcommand{\op}[1]{\buildrel\mbox{\tiny o}\over{#1}}
\newcommand{\db}{\prime\prime}
\newcommand{\bsl}{\backslash}
\newcommand{\lnrm}{\|\!|}
\newcommand{\rnrm}{|\!\|}
\newcommand{\lb}{\lbrack\!\lbrack}
\newcommand{\rb}{\rbrack\!\rbrack}
\newcommand\la{\langle}
\newcommand\ra{\rangle}
\newcommand{\ev}{\equiv}
\newcommand{\nev}{\not\equiv}
\newcommand{\nn}{\mathbb{N}}
\newcommand{\qq}{\mathbb{Q}}
\newcommand{\zz}{\mathbb{Z}}
\newcommand{\rr}{\mathbb{R}}
\newcommand{\rn}{\rr^N}
\newcommand{\cc}{\mathbb{C}}
\newcommand{\id}{\mathbb{I}}
\newcommand{\bo}{\mathbb{O}}

\newcommand{\amsb}[1]{\mathbb{#1}}
\newcommand{\mcl}[1]{\mathcal{#1}}
\newcommand{\bl}[1]{\mathbf{#1}}
\newcommand{\ov}[1]{\overline{#1}}
\newcommand{\wt}[1]{\widetilde{#1}}
\newcommand{\wh}[1]{\widehat{#1}}

\newcommand{\llra}{\leftrightarrow}
\newcommand{\lra}{\longrightarrow}
\newcommand{\LLR}{\Longleftrightarrow}
\newcommand{\LRA}{\Longrightarrow}
\newcommand{\LLA}{\Longleftarrow}


\newcommand{\bbox}{\vrule height.6em width.6em 
depth0em} 
\newcommand{\os}{\vbox{\hrule \hbox{\vrule 
height.6em depth0pt 
\hskip.6em \vrule height.6em depth0em}
\hrule}} 


\newcommand{\Ker}{\operatorname{Ker}}
\newcommand{\Imm}{\operatorname{Im}}
\newcommand{\rank}{\operatorname{rank}}
\newcommand{\dvg}{\operatorname{div}}
\newcommand{\curl}{\operatorname{curl}}
\newcommand{\supp}{\operatorname{supp}}
\newcommand{\essup}{\operatornamewithlimits{ess\,sup}}
\newcommand{\essinf}{\operatornamewithlimits{ess\,inf}}
\newcommand{\essosc}{\operatornamewithlimits{ess\,osc}}
\newcommand{\osc}{\operatornamewithlimits{osc}}
\newcommand{\sign}{\operatorname{sign}}
\newcommand{\loc}{\operatorname{loc}}
\newcommand{\diam}{\operatorname{diam}}
\newcommand{\dist}{\operatorname{dist}}
\newcommand{\card}{\operatorname{card}}
\newcommand{\meas}{\operatorname{meas}}
\newcommand{\spn}{\operatorname{span}}
\newcommand{\dtm}{\operatorname{det}}
%


\newcommand{\overlim}{\mathop{\overline{\lim}}\limits}
\newcommand{\underlim}{\mathop{\underline{\lim}}\limits}
\newcommand{\ttop}[2]{\genfrac{}{}{0pt}{}{#1}{#2}}
\newcommand{\bcu}{\mathop{\txty{\bigcup}}\limits}
\newcommand{\bca}{\mathop{\txty{\bigcap}}\limits}
\newcommand{\bsu}{\mathop{\txty{\sum}}\limits}
\newcommand{\pro}{\mathop{\txty{\prod}}\limits}


\newcommand{\pl}{\partial}
\newcommand{\ptt}{\frac{\pl}{\pl t}}
\newcommand{\ppx}{\frac\pl{\pl x}}
\newcommand{\dds}{\frac d{ds}}
\newcommand{\ddt}{\frac d{dt}}

\newcommand{\intl}{\int\limits}
\newcommand{\iintl}{\iint\limits}
\def\Xint#1{\mathchoice
    {\XXint\displaystyle\textstyle{#1}}%
    {\XXint\textstyle\scriptstyle{#1}}%
    {\XXint\scriptstyle\scriptscriptstyle{#1}}%
    {\XXint\scriptscriptstyle\scriptscriptstyle{#1}}%
    \!\int}
\def\XXint#1#2#3{\setbox0=\hbox{$#1{#2#3}{\int}$}
    \vcenter{\hbox{$#2#3$}}\kern-0.5\wd0}
\def\bint{\Xint-}
\def\dashint{\Xint{\raise4pt\hbox to7pt{\hrulefill}}}
\def\dashiint{\bint\kern-0.15cm\bint}

\newcommand{\ovl}[3]{\int_{#1}^{#2}\kern-#3pt\raise4pt\hbox to7pt{\hrulefill}\ }

\newcommand{\ovll}[3]{\intl_{#1}^{#2}\kern-#3pt\raise4pt\hbox to7pt{\hrulefill}\ }

\newcommand{\tvl}[2]{\iint_{#1}\kern-#2pt\raise4pt\hbox to7pt{\hrulefill}\ }



\newcommand{\omt}{\Om_T}
\newcommand{\plo}{\partial\Omega}
\newcommand{\ovo}{\bar{\Om} }

%
\newcommand{\ci}[1]{C^\infty\!\left({#1}\right)}
\newcommand{\cio}[1]{C_o^\infty\!\left({#1}\right)}
\newcommand{\lloc}[1]{L_{\loc}\!\left({#1}\right)}
\newcommand{\xy}{|x-y|}


\newcommand{\intom}{\intl_{\Om}}
\newcommand{\intbo}{\intl_{\plo}}
\newcommand{\inom}{\int_{\Om}}
\newcommand{\inbo}{\int_{\plo}}
\newcommand{\intrn}{\intl_{\rn}}


\newcommand{\bye}{\end{document}}



%
%
%

%
%
%
%
%

%
%

%% file: wiener_mac.tex
\def\ple{\partial E}
\def\om{\omega}
\def\Om{\Omega}
\def\lam{\lambda}
\def\phi{\varphi}
\def\extr{\Big\arrowvert}
\def\iint{\int\mskip-30mu\int}
\def\Iint{\int\mskip-30mu\int\limits}
\def\eps{\epsilon}
\def\varep{\epsilon}
\def\mum{\mu_-}
\def\mup{\mu_+}
\def\muduer{\mu(2R)}
\def\mur{\mu(R)}
\def\muq{\mu_{Q(R)}}
\def\umenkm{(u-k)_-}
\def\umenkp{(u-k)_+}
\def\umenpm{(u-k)_{\pm}}
\def\capparo{K_{\rho}}
\def\capsigro{K_{\sigma\rho}}
\def\niuno{\nu_1}
\def\nidue{\nu_2}
\def\nip{\nu^+}
\def\nim{\nu^-}
\def\csi{\xi}
\def\csip{\xi^+}
\def\csim{\xi^-}
\def\cspm{\xi^{\pm}}
\def\del{\delta}
\def\lam{\lambda}
\def\akrot{A_{k,\rho}}
\def\alrot{A_{l,\rho}}
\def\al{\alpha}
\def\w1p{W^{1,p}(\rn)}
\def\degpin{[DG]^+_p(E)}
\def\degmin{[DG]^-_p(E)}
\def\degpmin{[DG]^{\pm}_p(E)}
\def\degin{[DG]_p(E)}
\def\degpbo{[DG]^+_p(h;\bar E)}
\def\degmbo{[DG]^-_p(h;\bar E)}
\def\degpmbo{[DG]^{\pm}_p(h;\bar E)}
\def\degbo{[DG]_p(h;\bar E)}
\def\poten{{\cal U}_K}
\def\media{\mkern12mu\hbox{\vrule height4pt
           depth-3.2pt width5pt}\mkern-16.5mu\int
           \nolimits}
\def\lap{L^{\alpha,p}}
\def\wap{W^{\alpha,p}}
\def\capac{\hbox{cap}}
\newcommand{\dg}{[DG]_p^\pm(E)}
\newcommand{\gdg}{[GDG]_p^\pm(E)}
\newcommand{\dgp}{[DG]_p^+(E)}
\newcommand{\gdgp}{[GDG]_p^+(E)}
\newcommand{\dgm}{[DG]_p^-(E)}
\newcommand{\gdgm}{[GDG]_p^-(E)}
\newcommand{\hk}{\chi_{[k>h]}}